\documentclass{article}

\usepackage{PRIMEarxiv}

\usepackage[utf8]{inputenc} % allow utf-8 input
\usepackage[T1]{fontenc}    % use 8-bit T1 fonts
\usepackage{hyperref}       % hyperlinks
\usepackage{url}            % simple URL typesetting
\usepackage{booktabs}       % professional-quality tables
\usepackage{amsfonts}       % blackboard math symbols
\usepackage{nicefrac}       % compact symbols for 1/2, etc.
\usepackage{microtype}      % microtypography
\usepackage{lipsum}
\usepackage{fancyhdr}       % header
\usepackage{graphicx}       % graphics
\graphicspath{{media/}}     % organize your images and other figures under media/ folder

\usepackage{amsmath}
\usepackage{float}

\title{A NOTE ON AN OPEN CONJECTURE IN RATIONAL DYNAMICAL SYSTEMS
\thanks{\textit{\underline{Citation}}: 
\textbf{Zeraoulia Rafik and Alvaro humberto Salas, Open conjecture in rational dynamical system}
}}

\author{
  Zeraoulia Rafik \\
  University of batna2 (fesidis).Algeria \\
  Departement of mathematics\\
  yabous ,khenchela\\
  \texttt{\{Correspending Author\}r.zeraoulia@univ-batna2.dz} \\
   \And
  Alvaro Humberto Salas \\
  University national de Colombia,Bogota \\
  Department of Mathematics\\
  FIZMAKO Research Group \\
  \texttt{ahsalass@unal.edu.co}
 \\
}

\begin{document}
\maketitle

\begin{abstract}
Recently ,mathematicians have been interested in studying the theory of discrete dynamical system, specifically   difference equation, such that  considerable works about  discussing  the behavior properties  of its solutions (boundedness and unboundedness)  are discussed and published in many areas of mathematics which involves several interesting results and applications in applied mathematics and physics ,One of the most important discrete dynamics which is became of interest for researchers in the field is the rational dynamical system .In this paper we give a negative answer  to  the eight open conjecture in rational dynamical system  proposed by   G.Ladas and Palladino many years ago which states :
    
Assume $\alpha,\beta, \lambda \in [0,\infty)$. Then every positive solution of the difference equation \\: 

\begin{align*}
z_{n+1}=\frac{\alpha+z_{n}\beta +z_{n-1}\lambda}{z_{n-2}},\quad n=0,1,\ldots 
\end{align*}
is bounded if and only if $\beta=\lambda$. 
  We will use a construction of subenergy function and some properties of Todd's difference equation to disprove that conjecture in general.Some new results (Chebychev approximation) and analysis regarding that open conjecture are presented.

\end{abstract}

% keywords can be removed

\section{Introduction}
The theory of difference equations finds many applications
in almost all areas of natural science \cite{Dan:06}. 
increasingly clearly
emerges the fundamental role that difference equations with
discrete and continuous argument is played for understanding
nonlinear dynamics and phenomena also it is used for combinatorics and in the approximation of solutions of partial differential equations \cite{far:17}. The increased interest in difference equations is partly
due to their ease of handling. A minimum is enough
computing and graphical tools to see how the
solution of difference equations trace their bifurcations with
changing parameters \cite{Josef:08}. Thus opens a complex understanding as well invariant manifolds for linear and nonlinear dynamical systems.

\noindent  Let us define the sequence: $x=x_n$, $n\in \mathbb{Z^{+}} $  every  term
which is related to the previous recurrence relation
\begin{equation}\label{eq_R1}
    x_n=f(n,x_{n-1},x_{n-2},x_{n-3},\cdots x_{n-k})
\end{equation}
With a fixed $k >0$ ,the autonomous  variable $n$ changes
dis-continuously , and the formula  defined in  (\ref{eq_R1}) are called difference equations with
discrete arguments. If $x$ a is a function of continuous argument $m \in \mathbb{R^{+}}$ then the relation:
\begin{equation}\label{eq_R2}
    x(m)=f(t,x(m-1),x(m-2),\cdots x(m-k))
\end{equation}
is a difference equation with a continuous argument.
In practice, time usually  plays the role of an autonomous  variable, which
allows us to speak, respectively, of difference equations with
dis-continuous  and continuous time.Discrete time equations get up   when
the quantity $x$  under attention is recorded at some interval over
time. For example, if $x$ is the proportionate abundance (compactness,density)
of any biological kind , then as such an
interval, it is sometime  recomonded  to hold  the life time of one generation .
In sometimes  the relationship among $x_n$ and $x_{n-1}$  satisfactorily
is given by the first-order difference equation
\begin{equation}\label{eq_R3}
x_n =\lambda x_{n-1} (1-x_{n-1})
\end{equation}
(values $x_n$ as population density should not go out of
interval $[0,1]$, therefore, as it is easy to see, the parameter $x$ is,  the coefficient
reproduction should be enclosed in  $[0, 4]$)\cite{J.M:19}. If  $0<\lambda <1$ then the population is going  to zero at a rate
of power geometric sequence, if $1 < \lambda <  4$, then the comportement  of $x_n$ can
be  two of them  simple (perhapes stabilize over time or
become cyclic), and very complex (chaotic).
Difficulties in the behavior of $x_n$ arise due to the nonlinearity on the right side of 
equation (\ref{eq_R3}). it is not just a matter of nonlinearity, but 
that the right-hand side is essentially nonlinear, namely:the following  segment of  the increase $\bigg(\textit{interval}\quad (0,\frac12)\bigg)$ is followed by the segment
descending $\bigg(\textit{interval}\quad (\frac12,1)\bigg)$.
\noindent  It turns out that by studying real physical
problem, it is more convenient to first derive relations for finite
differences, make a passage to the limit, obtain differential
equations and only then by discretization in time and
space to arrive at difference schemes \cite{Ngoc:21}. Possibly due in part
precisely to these reasons, the development of the theory of difference equations, starting from
the end of the XVIII century, gradually lags behind the rapidly and multifaceted
 developing theory of differential equations as ordinary,
as well as in private derivatives.

\noindent  Attempts to understand the mechanisms of turbulence from  infinitesimal
point of view inevitably encounter various obstacles \cite{H:78},
caused by the need to solve the Navier-Stokes equations or
other nonlinear equations not inferior to them in complexity . Not
whether it is necessary to clarify the properties of turbulence completely different
equations that also reflect its discrete nature, which in
increasingly obvious lately? We mean such features of disturbance ,such as intermittency, the construction of various kinds,
consistent  structures, such as cyclones , etc. Recently formed
direction  structural turbulence \cite{Zambri:018}
confirms the adequacy of this opinion and  that suggestation  of the essence of the phenomenon. Nothing like
it is precisely the difference equations that must become such equations.
equations. However, from the properties of their solutions one can surprisingly clearly guess
many features of turbulence, and, as one would expect,
precisely those whose modeling with the help of differential
equations are the most difficult. For example, in the book \cite{Zambri:018} for modeling
sequences of operations for the formation of cyclones and vortices  of decreasing size,
the theory of hydrodynamic systems was developed. Each such
 system may include nonlinear ordinary differential
equation,  to obtain the subsequent scale of the vortices, 
the order has to be increased by three. As a result, the dimension of the system
is growing catastrophically and it becomes not at all easy to explore it.
simultaneously, the mechanism of the cascade operation itself can be implemented
already within the framework of only one dynamic   of the form
(\ref{eq_R3}).
\noindent
A  great  and interesting works were discussed and 
achieved by many researchers (published papers ,books ,notes,ect.. )\cite{Chri:095}  in the field of discrete rational dynamical systems such that it has been used to determine the behavior (  boundedness  ) of the solutions of a rational difference equations of the following form: 
\begin{equation}\label{eq}
x_{n+1}=\frac{\alpha+\beta x_{n}+\gamma x_{n-1}+\delta x_{n-2}}{A +B x_{n}+C x_{n-1}+D x_{n-1}}
\end{equation}
With nonnegative parameters $ \alpha ,\beta, \gamma ,\delta ,A,B,C$ and $D$  ,the main purpose  about studying the boundedness properties of the solutions of  the dynamics defined in (\ref{eq}) is to check  and to prove whether  solutions are still  bounded for all positive initial conditions or there exist some positive initial conditions where the solutions are  unbounded  .Dynamics of Third-Order Rational Difference Equations with Open Problems and Conjectures \cite{E. Camouzis:2008} deals with  large class of difference equations described by Equation (\ref{eq}),Some open problems related to (\ref{eq}) in which the boundedness properties were not known was recently solved in \cite{Peter:2019}, By the following assumption  $\delta=A=B=C=0$ with the variable change $x_n\to \dfrac{x_n}{D}$ , with $\alpha \geq 0,\beta >0,\gamma >0, D>0$ equation  (\ref{eq}) reduces to the following form :
\begin{equation}\label{Eq}
x_{n+1}=\frac{D\alpha+x_{n}\beta +x_{n-1}\gamma}{x_{n-2}},\quad n=0,1,\ldots , D\alpha=\alpha'
\end{equation}
It is shown in a paper by Lugo and Palladino \cite{Paladino:2009} that there exist unbounded solutions of (\ref{Eq})in the case that $0\leq \alpha <1$ and $0<\beta <\frac13$.Ying Sue Huang and Peter M. Knopf showed in \cite{Peter:2019} for $\alpha' \geq 0, \beta>0$ and if $\beta \neq 1$  there exist positive
initial conditions such that the solutions are unbounded except for
the case $\alpha'=0$ and $\beta >1$, In this paper we shall disprove the only if part of the eight conjecture of G.Ladas,Lugo and Palladino \cite{G.Ladas:12} such that we shall show using subenergy function and some numerical evidence (Mathematica Code) up to $10^{40}$ solutions that we have :
\begin{equation}\label{prof}
z_{n+1}=\frac{\alpha+z_{n}\beta +z_{n-1}\lambda}{z_{n-2}},\quad n=0,1,\ldots 
\end{equation}
implying
$\beta=\gamma$ 
is true ,but the converse is not .

\section{Conjecture}
Assume $\alpha,\beta, \lambda \in [0,\infty)$. Then every positive solution of the difference equation : 
\begin{equation}\label{eq_1}
z_{n+1}=\frac{\alpha+z_{n}\beta +z_{n-1}\lambda}{z_{n-2}},\quad n=0,1,\ldots 
\end{equation}
is bounded if and only if $\beta=\lambda$
\section{Proof}
Suppose that $\beta=\lambda>0$. Let $x_n:=z_n/\beta$ and $c:=\alpha/\beta^2$. Then the dynamics (\ref{eq_1}) can be rewritten as 

\begin{equation}\label{eq_2}
\qquad x_{n+1}=\frac{c +x_n+x_{n-1}}{x_{n-2}}
\end{equation}
(say for $n=2,3,\dots$), just with one parameter $c\ge0$ , the dynamic (\ref{eq_2})is exactly the Todd's difference equation ,in this case the equation is generally referred to by the cognomen “Todd’s equation” and possesses the invariant :
 
\begin{equation}\label{eq_3}
     \displaystyle(c +x_n+x_{n-1}+x_{n-2})\biggl(1+\dfrac{1}{x_n}\biggr)\biggl(1+\dfrac{1}{x_{n-1}}\biggr)\biggl(1+\dfrac{1}{x_{n-2}}\biggr)=\text{constant}
\end{equation}
The invariants of difference equations play an important role in understanding the stability and qualitative behavior of their solutions. To be more precise, if the invariant is a bounded manifold \cite{Y.KATO:2004},
then the solution is also bounded, Recently Hirota et al. \cite{R:2001}  found two conserved quantities $H_n^{1}$ and $H_n^{2}$ for the third- order Lyness equation , note that Lyness equation is a special case of equation (\ref{eq_2}) such that $c=1$ ,the two quantities are independents and One of the conserved quantities is the same form as that of (\ref{eq_3}) ,Both of two conserved quantities formula  were derived from   discretization of an anharmonic oscillator namely using its equation of its motion see the first equation here \cite{ R:2001}, we may consider those conserved quantities as a conserved subenergy of anharmonic oscillator, this means that (\ref{eq_3}) presents a sub energy function of that  anharmonic oscillator , To prove the "if" part of the conjecture  it would be enough to construct  for each nonnegative $c$, a "subenergy" function \cite{V.Z:2015} $f_c\colon(0,\infty)^3\to\mathbb{R}$ such that :
\begin{equation}\label{eq_4}
\qquad f_c (x_0,x_1,x_2)\to\infty\quad\text{as}\quad x_0+x_1+x_2\to\infty
\end{equation}
Note that the subenergy function is the invariant of the third difference equation ,namely, the dynamics (\ref{eq_2}) , if we assume that :
\begin{equation}\label{eq_5}
f_c (x_n,x_{n-1},x_{n-2})=\displaystyle(c +x_n+x_{n-1}+x_{n-2})\biggl(1+\dfrac{1}{x_n}\biggr)\biggl(1+\dfrac{1}{x_{n-1}}\biggr)\biggl(1+\dfrac{1}{x_{n-2}}\biggr)=\text{constant}n\geq0
\end{equation}

 then the condition (\ref{eq_4}) is satisfied in (\ref{eq_5}) .see Lemma2 in (\cite{M.ROG:2009},p.4) .For RHS of (\ref{eq_5}) see also \textbf{Theorem2.1} in (\cite{ M.ROG:2009},p.31) ,And  Since the invariant of the dynamic of (\ref{eq_2}) is constant then $f_c$  could be referred to as the conservation of energy along the path of the dynamical system.For some natural $k$ and all $x=(x_0,x_1,x_2)\in(0,\infty)^3$ one has the "subenergy" inequality
$f_c(T^k x)\le f_c(x)$, where $Tx:=(x_1,x_2,x_3)$, with $x_3=\frac{c+x_2+x_1}{x_0}$, according to the dynamics. Of course, $T^k$ is the $k$th power of the operator $T$. For $k=1$, the sub-energy inequality is the functional inequality 

\begin{equation}\label{eq_6}
\qquad f_c\Big(x_1,x_2,\frac{c+x_2+x_1}{x_0}\Big)\le f_c(x_0,x_1,x_2) \quad
\text{for all positive $x_0,x_1,x_2$, }
\end{equation}
To construct a subenergy function, one might want to start with some easy function $f_{c,0}$ such that $f_{c,0}(x_0,x_1,x_2)\to\infty$ as $x_0+x_1+x_2\to\infty$, and then consider something like $f_{c,0}\vee(f_{c,0}\circ T^k)\vee(f_{c,0}\circ T^{2k})\vee\dots$,Inequality (\ref{eq_6}) which can be clearly restated in the following more symmetric form: 
\begin{equation}
\qquad 
x_0x_3=c+x_1+x_2\implies f_c(x_1,x_2,x_3)\le f_c(x_0,x_1,x_2)  
\end{equation}
for all positive real $x_0,x_1,x_2,x_3$.  condition $x_0+x_1+x_2\to\infty$ in (\ref{eq_4}) can be replaced by any one of the following (stronger) conditions: (i) $x_0\to\infty$ or (ii) $x_1\to\infty$ or (iii) $x_2\to\infty$; this of course will replace condition (\ref{eq_4}) by a weaker condition, which will make it easier to construct a sub-energy function $f_c$ ,Here are details: Suppose that (\ref{eq_6}) holds for some function $f_c$ such that $f_c(x_0,x_1,x_2)\to\infty$ as $x_0\to\infty$. Suppose that, nonetheless, a positive sequence $(x_0,x_1,\dots)$ satisfying condition (\ref{eq_2}) is unbounded, so that, as $k\to\infty$, one has $x_{n_k}\to\infty$ for some sequence $(n_k)$ of natural numbers. Then $f_c(x_{n_k},x_{1+n_k},x_{2+n_k})\to\infty$ as $k\to\infty$. This contradicts (\ref{eq_4}), which implies, by induction, that $f_c(x_n,x_{1+n},x_{2+n})\le f_c(x_0,x_1,x_2)$ for all natural $n$. Quite similarly ,one can do with (ii) $x_1\to\infty$ or (iii) $x_2\to\infty$ in place of (i) $x_0\to\infty$. 

Moreover , instead of the dynamics of the triples $(x_n,x_{1+n},x_{2+n})$ one can consider the corresponding dynamics (in $n$) of the consecutive $m$-tuples $(x_n,\dots,x_{m-1+n})$ for any fixed natural $m$. 

Moreover , instead of inequality $f_c(x_1,x_2,x_3)\le f_c(x_0,x_1,x_2)$ in (\ref{eq_4}), one may consider a weaker inequality like $f_c(x_2,x_3,x_4)\le f_c(x_0,x_1,x_2)\vee f_c(x_1,x_2,x_3)$ for all positive $x_0,\dots,x_4$ satisfying condition (\ref{eq_2}), thanks to the invariant of Todd's difference equation (\ref{eq_5}) which it is defined in our case to be  a subenergy function such that it is easy to see that the if part of the conjecture would follow since the subenergy $f_c$ is always found .In   (\cite{E.A:2002}, p.35) Authors showed that every positive solution of dynamics (\ref{eq_2}) with invariant are bounded and persist this result is the affirmation that invariant must be a constant subenergy function which it is always found for all positive initial conditions \cite{Anna:2006}
 One can try to do the "only if" part in a similar manner. Suppose that $0<\beta\ne\lambda>0$. Let $u_n:=z_n/\sqrt{\beta\lambda}$, $c:=\alpha/(\beta\lambda)$, and $a:=\sqrt{\beta/\lambda}\ne1$. Then the dynamics (\ref{eq_2}) can be rewritten as: 
\begin{equation}\label{eq_7}
\qquad u_{n+1}=\frac{c+au_n+u_{n-1}/a}{u_{n-2}}, 
\end{equation}
Just with two parameters, $c \ge0$ and $a>0$. Suppose one can construct, for each pair $(c,a)\in[0,\infty)\times\big((0,\infty)\setminus\{1\}\big)$ and some $\rho=\rho_{c,a}\in(1,\infty)$, a "$\rho$-super-energy" function $g=g_{a,c;\rho}\colon(0,\infty)^3\to(0,\infty)$ such that $g$ is bounded on each bounded subset of $(0,\infty)^3$ and 

\begin{equation}\label{eq_8}
\qquad g\Big(u_1,u_2,\frac{c+au_2+u_1/a}{u_0}\Big)\ge\rho\, g(u_0,u_1,u_2)\quad 
\text{for all positive $u_0,u_1,u_2$.}
\end{equation}
Then, by induction, $g(u_n,u_{1+n},u_{2+n})\ge\rho^n g(u_0,u_1,u_2)\to\infty$ as $n\to\infty$, for any sequence $(u_n)$ satisfying (\ref{eq_7}). Therefore, and because $g$ is bounded on each bounded subset of $(0,\infty)^3$, it would follow that the sequence $(u_n)$ is unbounded. 

For any pair $(c,a)\in[0,\infty)\times(0,\infty))$ and any $\rho\in(1,\infty)$, there is no  "$\rho$-super-energy" function $g\colon(0,\infty)^3\to(0,\infty)$. This follows because the point $(u_{a,c},u_{a,c},u_{a,c})$ with 
$u_{a,c}:=\dfrac{1+a^2+\sqrt{a^4+4 a^2 c+2 a^2+1}}{2 a}$ is a fixed point (in fact, the only fixed point) of the map $T$ given by the formula $T(u_0,u_1,u_2)=\Big(u_1,u_2,\dfrac{c+au_2+u_1/a}{u_0}\Big)$. (If $a\ne1$, then this point is the only fixed point \cite{Esmmir:2019} of the map $T^2$ as well.)

This also disproves, in general, the "only if" part of the conjecture defined in (\ref{eq_1})

However, One may now try to amend this conjecture by excluding the initial point $(u_{a,c},u_{a,c},u_{a,c})$. Then, accordingly, the definition of a "$\rho$-super-energy" function would have it defined on a subset (say $S$) of the set $(0,\infty)^3\setminus\{(u_{a,c},u_{a,c},u_{a,c})\}$, instead of $(0,\infty)^3$; such a subset may be allowed to depend on the choice of the initial point $(u_0,u_1,u_2)$, say on its distance from the fixed point $(u_{a,c},u_{a,c},u_{a,c})$, and one would then have to also prove that $S$ is invariant under the map $T$ .

\section{Analysis and discussion}

\textbf{Case1:}

\begin{figure}[H]
    \centering
    \includegraphics[width=0.6 \textwidth]{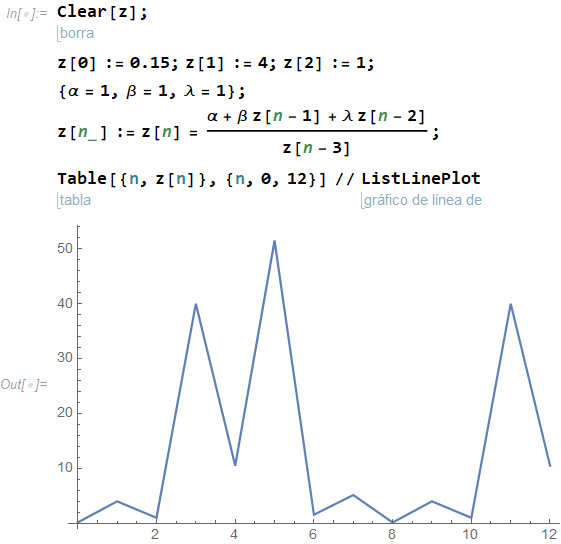}
    \caption{Bounded solutions for Todd dynamics in the case  $\alpha=\beta=\gamma$}
    \label{fig:my_1}
\end{figure}

Let us  try $\alpha=\beta=\gamma$  are the  positive initial conditions  in the dynamic (\ref{prof}),  dynamics  becomes a well- known third-order difference equation ,namely ,Todd difference equation,and all its solutions are bounded in that case ,see Figure \ref{fig:my_1} .

\begin{figure}[H]
    \centering
    \includegraphics[width=0.6\textwidth]{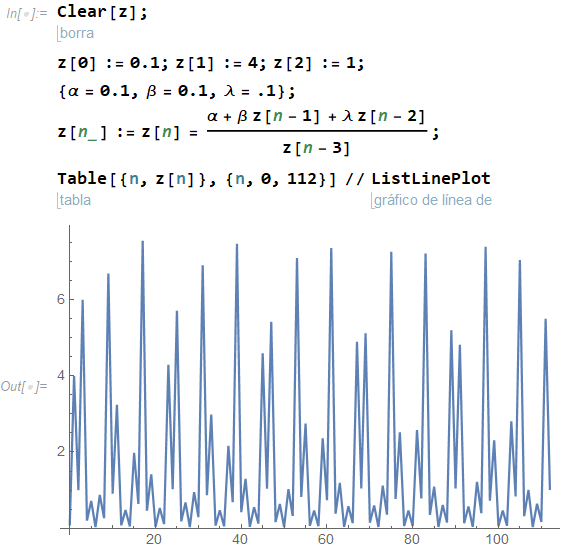}
    \caption{Bounded solution for the dynamic (\ref{prof}) in the case $\alpha=\beta <\gamma=1$}
    \label{fig:my_2}
\end{figure}

\textbf{case2}

let us  try now :$\alpha=\beta <\gamma=1$
In this case we still have a bounded solution as shown in Figure \ref{fig:my_2}

\begin{figure}[H]
    \centering
    \includegraphics[width=0.6\textwidth]{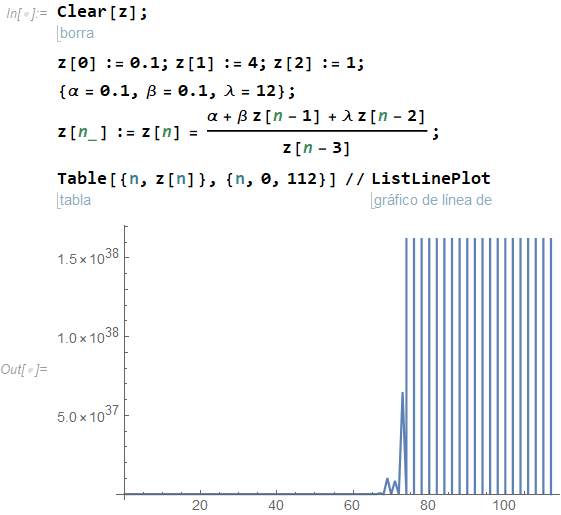}
    \caption{Bounded solution for the dynamic (\ref{prof}) in the case $\alpha=\beta <\gamma=1$}
    \label{fig:my_labe3}
\end{figure}

\textbf{case3}

for this case we may try :$\alpha=\beta <\gamma,\gamma >1$ , we get unbounded solutions as shown in figure \ref{fig:my_labe3}

This case indicate a lot about the only if part of the conjecture such that  one of the noted indications is that we should have  also  $\sigma$ to be lie  in the range $(0,1]$ for the conjecture to be hold  ,and one wishes to try  more values for $\alpha$ and $\beta$ .

\section{Chybeshev Approximation for bounds}

\noindent We have noted an interesting behavior of solutions of the dynamics defined in  (\ref{prof}) such that bounds behave linearly using Chebychev approximation as shown in Figure

\begin{figure}[H]
    \centering
    \includegraphics[width=0.6\textwidth]{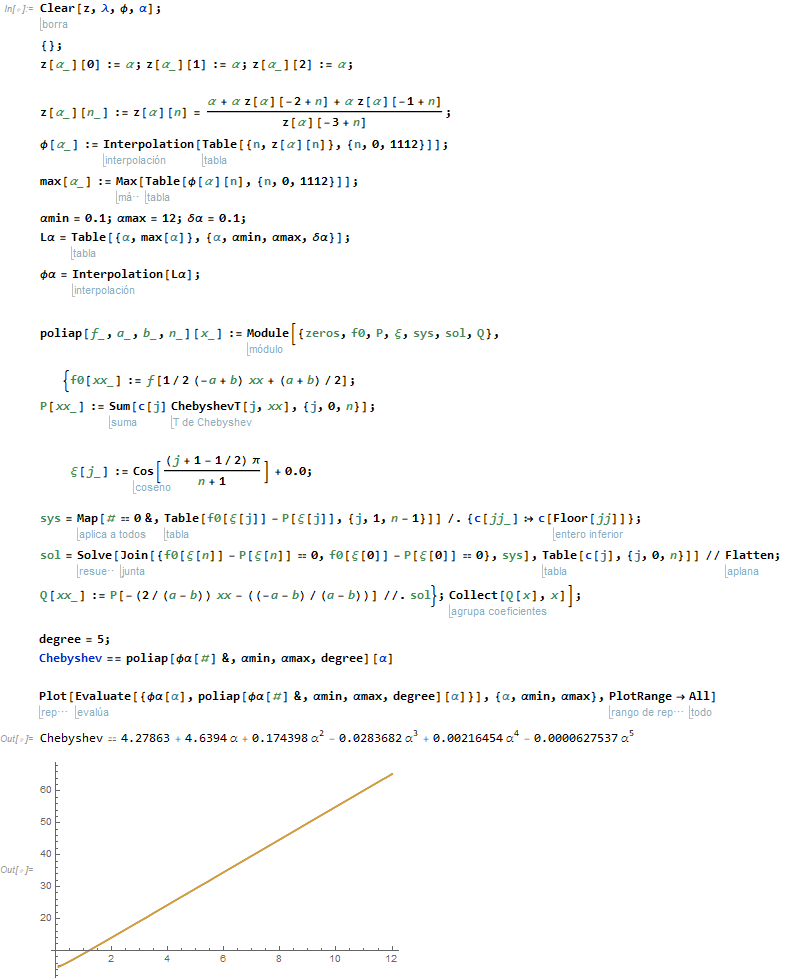}
    \caption{linear behavior of bounded solutions of the  dynamics (\ref{prof})using Chebychev approximation}
    \label{fig:my_4}
\end{figure}

\section{Conclusion}

Usually Lyaponov theory and some advanced theories in differential equations are unable to give an affirmative answer to an arbitrary challenged problem in rational dynamical system,  
whereas we may find a ready  affirmative answer to the challenged problem in physics using some interpretations which are depending  on the behavior and  properties of high energy functions like super energy function and hamiltonian operator .
\bigskip

\textbf{Data Availability
}
No data were used to support this paper.

\bigskip

\textbf{Conflicts of Interest
}
The authors declare no conflicts of interest with regard to any individual or organization.

\bigskip

\textbf{Acknowledgments}
Thanks to Allah who helped me to do that modest work ,This work was supported by the funding from the employer Universidad Nacional de Colombia and university of Batna2 .Algeria,I would like to express our deep gratitude to Professors Iosif Pinelis and My Co.Author my best friend Alvaro h Salas  for improving the quality of that paper and  to my parents(Abd alkarim and djemila) and all my great salutations to my lovely wife and to my second heart my son Taha Abd Aldjalil and all my brothers(saddam and khireddin and Djawad) and my sisters ,In particulary Arras Ta3i.


%% References


\bibliographystyle{elsarticle-harv}
%%\bibliography{riaibib}



%%\begin{thebibliography}{01}



%%\bibitem[{Able(1945)}]{Abl:45}
%%Able, B., 1945. Nombre del artículo. Nombre de la revista 35, 123--126.

\end{document}